\documentclass[12pt,a4paper]{article}
\usepackage{epsfig,latexsym,amsfonts,amssymb,amsmath,amscd,
graphics,theorem} 
\theoremstyle{plain}
\newtheorem{lemma}{Lemma}[section]
\newtheorem{proposition}[lemma]{Proposition}
\newtheorem{remark}[lemma]{Remark}

\newtheorem{theorem}[lemma]{Theorem}

{\theorembodyfont{\rmfamily} 
\font\ncsc=cmcsc10  \font\ntt=cmtt12
\begin{document}
\baselineskip=18pt
\newcommand{\pperp}{\hbox{$\perp\hskip-6pt\perp$}}
\newcommand{\ssim}{\hbox{$\hskip-2pt\sim$}}
\newcommand{\N}{{\mathbb N}}
\newcommand{\A}{{\mathbb A}}
\newcommand{\Z}{{\mathbb Z}}
\newcommand{\R}{{\mathbb R}}
\newcommand{\C}{{\mathbb C}}
\newcommand{\Q}{{\mathbb Q}}
\newcommand{\PP}{{\mathbb P}}
\newcommand{\mnote}{\marginpar}
\newcommand{\Id}{{\operatorname{Id}}}
\newcommand{\oeps}{{\overline\eps}}
\newcommand{\oDel}{{\widetilde\Del}}
\newcommand{\real}{{\operatorname{Re}}}
\newcommand{\conv}{{\operatorname{conv}}}
\newcommand{\Span}{{\operatorname{Span}}}
\newcommand{\Ker}{{\operatorname{Ker}}}
\newcommand{\Hyp}{{\operatorname{Hyp}}}
\newcommand{\Fix}{{\operatorname{Fix}}}
\newcommand{\sign}{{\operatorname{sign}}}
\newcommand{\Tors}{{\operatorname{Tors}}}
\newcommand{\oi}{{\overline i}}
\newcommand{\oj}{{\overline j}}
\newcommand{\ob}{{\overline b}}
\newcommand{\os}{{\overline s}}
\newcommand{\oa}{{\overline a}}
\newcommand{\oy}{{\overline y}}
\newcommand{\ow}{{\overline w}}
\newcommand{\ou}{{\overline u}}
\newcommand{\ot}{{\overline t}}
\newcommand{\oz}{{\overline z}}
\newcommand{\bw}{{\boldsymbol w}}
\newcommand{\bx}{{\boldsymbol x}}
\newcommand{\bu}{{\boldsymbol u}}
\newcommand{\bz}{{\boldsymbol z}}
\newcommand{\eps}{{\varepsilon}}
\newcommand{\proofend}{\hfill$\Box$\bigskip}
\newcommand{\Int}{{\operatorname{Int}}}
\newcommand{\pr}{{\operatorname{Pr}}}
\newcommand{\grad}{{\operatorname{grad}}}
\newcommand{\rk}{{\operatorname{rk}}}
\newcommand{\im}{{\operatorname{Im}}}
\newcommand{\sk}{{\operatorname{sk}}}
\newcommand{\const}{{\operatorname{const}}}
\newcommand{\Sing}{{\operatorname{Sing}}}
\newcommand{\conj}{{\operatorname{Conj}}}
\newcommand{\Pic}{{\operatorname{Pic}}}
\newcommand{\Crit}{{\operatorname{Crit}}}
\newcommand{\Ch}{{\operatorname{Ch}}}
\newcommand{\discr}{{\operatorname{discr}}}
\newcommand{\Tor}{{\operatorname{Tor}}}
\newcommand{\Conj}{{\operatorname{Conj}}}
\newcommand{\val}{{\operatorname{val}}}
\newcommand{\Val}{{\operatorname{Val}}}
\newcommand{\defect}{{\operatorname{def}}}
\newcommand{\tmu}{{\C\mu}}
\newcommand{\ov}{{\overline v}}
\newcommand{\ox}{{\overline{x}}}
\newcommand{\tet}{{\theta}}
\newcommand{\Del}{{\Delta}}
\newcommand{\bet}{{\beta}}
\newcommand{\kap}{{\kappa}}
\newcommand{\del}{{\delta}}
\newcommand{\sig}{{\sigma}}
\newcommand{\alp}{{\alpha}}
\newcommand{\Sig}{{\Sigma}}
\newcommand{\Gam}{{\Gamma}}
\newcommand{\gam}{{\gamma}}
\newcommand{\Lam}{{\Lambda}}
\newcommand{\lam}{{\lambda}}
\newcommand{\SC}{{SC}}
\newcommand{\MC}{{MC}}
\newcommand{\nek}{{,...,}}
\newcommand{\cim}{{c_{\mbox{\rm im}}}}
\newcommand{\mathto}{\mathop{\to}}

\newcommand{\w}{{\omega}}

\title{Welschinger invariant and enumeration of\\
real rational curves}
\author{Ilia Itenberg \and Viatcheslav Kharlamov
\and Eugenii Shustin}
\date{}
\maketitle
\begin{abstract}
Welschinger's invariant bounds from below the number of real
rational curves through a given generic collection of real points
in the real projective plane. We estimate this invariant using
Mikhalkin's approach which deals with a corresponding count of
tropical curves. In particular, our estimate implies that, for any
positive integer $d$, there exists a real rational curve of degree
$d$ through any collection of $3d-1$ real points in the projective
plane, and, moreover, asymptotically in the logarithmic scale at
least one third of the complex plane rational curves through a
generic point collection are real. We also obtain similar results
for curves on other toric Del Pezzo surfaces.
\end{abstract}

\begin{figure}[htb]
\hfill\includegraphics[height=2.5cm,width=8cm,angle=0,draft=false]{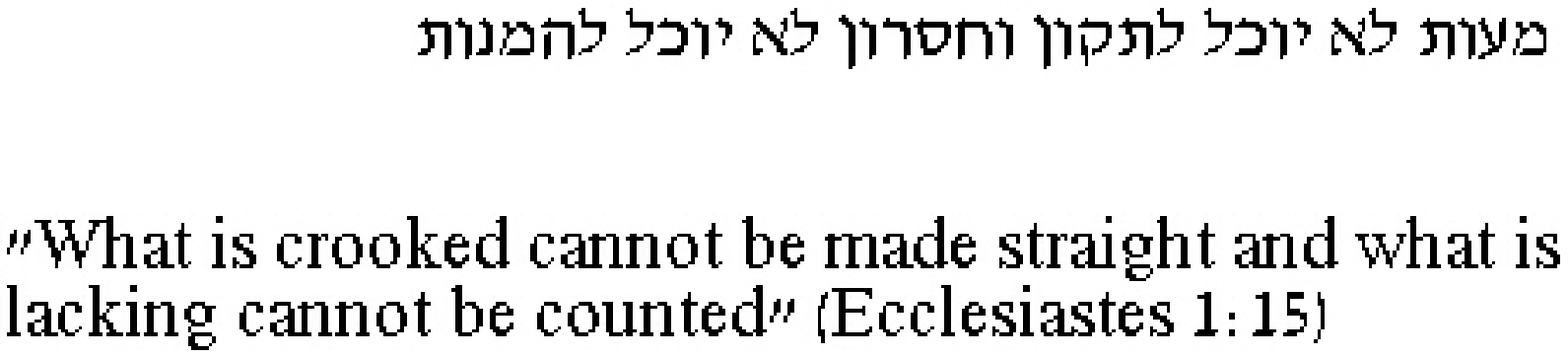}
\end{figure}

\section{Introduction}\label{intro}

In contrast to spectacular achievements in complex enumeration
geometry, real enumerative geometry remained in its almost
embryonic state, with Sottile's~\cite{So} and
Eremenko-Gabrielov's~\cite{EG} results on the real Schubert
calculus as the only serious exception. One of the difficulties in
real enumerative geometry is the lack of invariants: the number of
real objects usually varies along the parameter space. One of the
simplest known examples \cite{DK} showing this phenomenon is the
following: while over $\C$ the number of rational plane cubics
through given $8$ generic points is $12$, the number of real ones
can take three different values, $8$, $10$, and $12$, depending on
the positions of the $8$ real points. One of the main questions in
real enumerative geometry is that of the lower bounds on the
number of real solutions (an upper bound provided by the number of
complex solutions gives rise to another interesting question, that
is of the sharpness of such a bound, see the discussion of "total
reality" in~\cite{So} and further references therein).

In the present note we focus on the enumeration of real rational
curves passing through fixed real points on an algebraic surface,
a real counterpart of the complex Gromov-Witten theory, and, in
particular, on the following question: {\it whether through any
generic $3d-1$ points in the real plane always passes a real
rational curve of degree $d$\rm }? (The respective number of
complex rational curves~\cite{KM} is even for every $d \ge 3$, so
the existence of required real curves does not immediately follow
from the computation in the complex case.)

Till recently, the answer to the above question was not known even
for degree~$4$. The situation has radically changed after the
discovery by J.-Y.~Welschinger~\cite{W, W1} of a way of attributing
weights to real solutions making the number of solutions counted
with the weights to be independent of the configuration of
points. Since the weights take values $\pm 1$ exclusively,
the absolute value of Welschinger's invariant $W_d$
immediately provides
a lower bound on the number $R_d$ of real solutions:
$R_d\ge|W_d|$. The next natural question arises: {\it how
non-trivial is Welschinger's bound\rm }?

In the present note we show that
Welschinger's
bound
implies the following statement.
\begin{theorem}\label{t3}$\;$
For any
integer $d\ge 1$, through any $3d-1$ generic points in $\R P^2$
there can be traced at least $d!/2$ real rational curves of degree
$d$.
\end{theorem}

As a corollary, we obtain an affirmative answer to the
aforementioned question:
there always
exists
at least one real rational curve going through
the given real
points
(such
an existence statement extends from generic point data to an
arbitrary one at least if reducible rational curves are taken into
consideration as well).

To establish the non-triviality of the Welschinger invariant and
to estimate it, we apply another key ingredient, Mikhalkin's
approach to counting nodal curves passing through specific
configurations of points \cite{Mi}. We also point out that the
Welschinger invariant does not generalize directly even to the
case of elliptic curves (section \ref{sec4}), thereby leaving open
the enumeration problem of real plane curves of positive genera.

Note that the
expression given in Theorem \ref{t3}
is not
the exact value
of
the Welschinger invariant. Using the same
methods one can perform explicit calculations for, say,  small
$d$. For example, in this way one gets $W_4=240$ and $W_5=18264$.

\begin{remark}
Comparing with the formulae in {\rm\cite{KM}}, one can see that
the logarithm of the lower bound in
Theorem~\ref{t3}
is
asymptotically equal to $\frac{1}{3}\log N^\C_{\PP^2,d}$, where
$N^\C_{\PP^2,d}$ is the number of complex rational plane curves of
degree~$d$ passing through $3d - 1$ generic points. {\rm (In fact,
$\log N^\C_{\PP^2,d}\sim 3d\log d$, as follows from the
inequalities $(3d-4)!\cdot 54^{-d}\le N^\C_{\PP^2,d}\le (3d-5)!$
which, in turn, easily follow from Kontsevich's recurrence
formula~{\rm\cite{KM}}.)}
\end{remark}

\vskip10pt

We treat in a similar way enumerative problems for real rational
curves on
other toric Del Pezzo surfaces
equipped with
their usual real structures,
and obtain the following statement
as a corollary of
Welschinger's bound.

Denote by~$Q$ the hyperboloid $\C P^1 \times \C P^1$, and by $P_k$,
$k = 1, 2, 3$ the complex projective plane with $k$ blown
up
generic real points. For
$P_k$, let~$L$ be the pull back of a generic straight line, and
$E_1,...,E_k$ the exceptional divisors.

\begin{theorem}\label{t3n} $\;$
Let $\Sigma$ be $Q$, $P_1$, $P_2$ or $P_3$,
and $D$
a very
ample divisor on
$\Sigma$. Then through any $c_1(\Sigma) \cdot D - 1$ generic real
points in $\Sigma$ there can be traced at least $\varrho$ real
rational curves belonging to the linear system $|D|$,
where
$$\varrho =
\begin{cases}
\max\{d_1, d_2\}, \quad & \text{\rm if} \ \Sigma = Q,
\hskip0.25cm D \ \text{\rm is of bi-degree} \ (d_1, d_2), \\
\frac{1}{2}\ d \ !/d_1 !,
\quad & \text{\rm if} \ \Sigma = P_1,
\ D \sim dL - d_1E_1, \ d > d_1>0,\\
(d-d_2)!/d_1!, \quad & \text{\rm if} \ \Sigma = P_2,
\ D \sim dL - d_1E_1 - d_2E_2, \\
& \hskip2.09cm d_1 + d_2 < d, \ \text{\rm and} \ d_1 \geq d_2>0, \\
(d-d_2-d_3)!/d_1!, \quad & \text{\rm if} \ \Sigma = P_3,
\ D \sim dL - d_1E_1 - d_2E_2 - d_3E_3, \\
& \hskip2.09cm d_1 + d_2 + d_3 \leq d, \ \text{\rm and} \
d_1 \geq d_2 \geq d_3>0, \\
(d_1!)/(d-d_2-d_3)! \quad & \text{\rm if} \ \Sigma = P_3,
\ D \sim dL - d_1E_1 - d_2E_2 - d_3E_3, \\
& \hskip2.09cm d_1 + d_2 + d_3 > d, \ \text{\rm and} \
d_1 \geq d_2 \geq d_3>0.
\end{cases}$$
\end{theorem}

In fact,
the Welschinger invariant counts real $J$-holomorphic curves
on real symplectic rational surfaces equipped with a generic almost complex
structure.
A literal use of this invariant for
an
integrable
structure on a surface is possible
if
for this surface
there is no difference between Gromov-Witten invariants and the count
of irreducible rational curves.
The latter condition is verified
for
Del Pezzo surfaces.
It is no longer the case for
a surface containing an exceptional divisor
$E$ with $E^2 \leq -2$.
On the other hand, Mikhalkin's approach
to counting nodal curves
applies only to toric surfaces.
That is why we restrict ourselves
to the case of toric Del Pezzo surfaces
$\C P^2$, $Q$ and $P_k$, $k = 1, 2, 3$.
The next natural question
would be on existence of a real rational curve
passing through a generic collection of real points
on a nonsingular toric surface, and in particular,
on a geometrically ruled surface $\Sigma_n$ with $n \geq 2$.
However, to solve this question one needs additional tools.

\vskip10pt

Due to the conformal uniqueness of symplectic structures on $\C
P^2$ and certain deformation invariance of Welschinger numbers
(see \cite{W, W1}),
Theorem \ref{t3}
extends to the corresponding
symplectic settings. In particular, by means of such a
generalization of
Theorem \ref{t3}
one can extend the
Hilbert-type inequalities for real plane algebraic curves (see,
for example, \cite{DK}) to the case of real pseudo-holomorphic
curves on the real symplectic projective plane. Define a partial
order on the set of connected components of a smooth curve in $\R
P^2$ in the following way: if a component $C$ is contained in the
disc bounded by a component $C'$, then $C'$ {\it dominates} $C$.

\begin{proposition}\label{t4}
For any nonsingular real pseudo-holomorphic curve of degree $d$ in
the real symplectic projective plane and any integer $s \ge 1$,
the total length of $3s-1$ disjoint linearly ordered chains of
real components of the curve does not exceed $ds/2$ if no element
of one chain dominates all the elements the other chains.
\end{proposition}

Note that in the integrable case
there is a proof
of
such a statement
(see \cite{DK}) which
is based on the possibility to
trace a connected real cubic (of genus $0$) through $8$ points and
a connected real quartic (of genus $3$) through $13$ points. A
similar proof could be proposed in the symplectic category. For
that, however, in addition to tracing a real rational
pseudo-holomorphic curve of degree $3$ through $8$ points (one of
particular Welschinger's results, \cite{W, W1}), we should know that
for a generic almost-complex structure the (one-dimensional)
family of pseudo-holomorphic curves of degree $4$ passing through
generic $13$ points contains an odd number of singular elements,
which seems to be still an open question.

\section{Tropical calculation of the Welschinger
invariant}\label{sec2}

\begin{figure}
\begin{center}
\epsfxsize 145mm \epsfbox{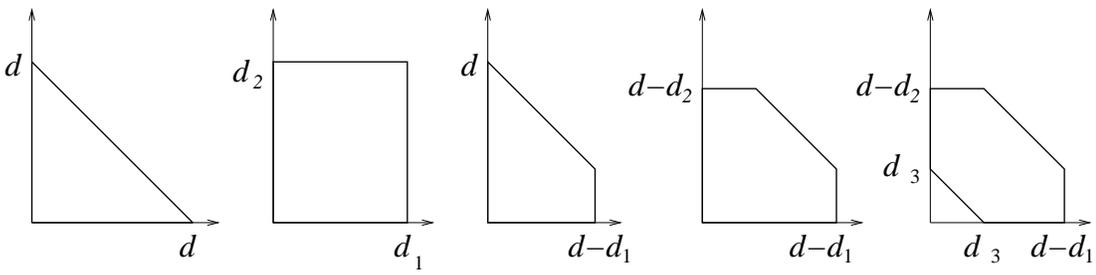}
\end{center}
\caption{Polygons associated with $\C P^2$, $Q$, $P_1$, $P_2$, and
$P_3$} \label{f2}
\end{figure}

Let $\Sig$ be $\C P^2$, $Q$, $P_1$, $P_2$ or $P_3$ equipped with
its standard real structure, and $D$ a very ample divisor on
$\Sigma$. The linear system $|D|$ is generated, with respect to
suitable affine coordinates, by monomials $x^iy^j$, where $(i,j)$
ranges over all the integer points ({\it i.e.}, points having
integer coordinates) of a convex lattice polygon~$\Delta$ of the
following form. If $\Sig= \C P^2$ and $D \sim d[\C P^1]$, then
$\Delta$ is the triangle with vertices $(0,0)$, $(d,0)$, $(0,d)$.
If $\Sig=Q$ and $D$ is of bi-degree $(d_1, d_2)$, then $\Delta$ is
the rectangle with vertices $(0,0)$, $(d_1,0)$, $(d_1,d_2)$,
$(0,d_2)$. If $\Sig=P_k$, $k=1,2,3$, and $D\sim dL-\sum_{i=1}^k
d_iE_i$, then $\Del$ is respectively the quadrangle with vertices
$(0,0)$, $(d-d_1,0)$, $(d-d_1,d_1)$, $(0,d)$, or the pentagon with
vertices $(0,0)$, $(d-d_1,0)$, $(d-d_1,d_1)$, $(d_2,d-d_2)$,
$(0,d-d_2)$, or the hexagon with vertices $(d_3,0)$, $(d-d_1,0)$,
$(d-d_1,d_1)$, $(d_2,d-d_2)$, $(0,d-d_2)$, $(0,d_3)$ (see Figure
\ref{f2}). Let $r(\Delta)$ be the number of integer points on the
boundary of $\Delta$ diminished by $1$, and $\delta(\Delta)$ be
the number of interior integer points of~$\Delta$. Note that
$r(\Delta) = c_1(\Sig)\cdot D-1$ and $\delta(\Delta)$ is the genus
of nonsingular representatives of $|D|$. As well known, the number
of curves of genus $0\le g\le\delta(\Delta)$ in $|D|$ passing
through $c_1(\Sig)\cdot D-1+g$ generic points is finite.

\medskip

{\bf A. Welschinger numbers.}\label{Welsch-inv} For any integer $0
\leq g \leq \delta(\Delta)$ and any set $\bw$ of $c_1(\Sig)\cdot
D-1+g$ generic distinct real points in $\Sig$, let
$N^\R_{\Sig,D,g}(\bw)$ be the number of irreducible real nodal
curves of genus $g$ in $|D|$ passing through all the points of
$w$, and let $N_{\Sig,D,g}^{even}(\bw)$ (resp.,
$N_{\Sig,D,g}^{odd}(\bw)$) be the number of real irreducible nodal
curves of genus $g$ in $|D|$ passing through all the points of
$\bw$ and having even (resp., odd) number of solitary nodes ({\it
i.e.}, double points locally given by $x^2+y^2=0$). Define the
{\it Welschinger number} as
$W_{\Sig,D,g}(\bw)=N_{\Sig,D,g}^{even}(\bw)-N_{\Sig,D,g}^{odd}(\bw)$.

\begin{theorem}\label{Welsch-theorem}{\rm (J.-Y.~Welschinger,
see~\cite{W, W1})}.
If $g=0$, then $W_{\Sig,D,g}(\bw)$
does not depend on the choice of the {\rm(}generic{\rm)}
set~$\bw$.
\end{theorem}

We call $W_{\Sig,D,0}(\bw)$ the {\it Welschinger invariant} and
denote it by $W_{\Sig,D}$. Clearly, $N^\R_{\Sig,D,0}(\bw)\ge
|W_{\Sig,D}|$. This Welschinger bound together with an explicit
estimate for $W_{\Sig,D}$, obtained below, implies
Theorems~\ref{t3}, \ref{t3n}. To get the estimate we use the
methods of tropical algebraic geometry \cite{Mi,Sh}.

{\bf B. Correspondence theorems.} Let $A$ be a finite collection
of integer points in $\R^2$, and $\nu: A \to \R$ a function.
Consider the function $\hat \nu: \R^2 \to \R$ defined by $\hat \nu
(x,y) = \max_{(i,j)\in A}\{ix + jy - \nu(i,j)\}$. The function
$\hat \nu$ is called the {\it Legendre transform} of~$\nu$. Note
that $\hat \nu$ is a piecewise-linear convex function. Consider
the corner locus $\Pi \subset \R^2$ of $\hat \nu$, {\it i.e.}, the
set where $\hat \nu$ is not smooth. The set $\Pi$ has a natural
structure of one-dimensional complex, and any edge of $\Pi$ is a
segment or a ray. Denote by $\Delta(A)$ the convex hull of~$A$.
Any connected component of the complement of $\Pi$ corresponds to
an integer point of $\Delta(A)$. Thus, the edges of $\Pi$ are
naturally equipped with positive integer weights. Namely, the
weight of an edge separating two connected components of the
complement of $\Pi$ is the (integer) length of the segment joining
the corresponding two integer points of $\Delta(A)$. The resulting
weighted complex is called the {\it tropical curve} associated
with the pair $(A, \nu)$; cf.~\cite{Mi2} and~\cite{St}. The
polygon $\Delta(A)$ is called a {\it Newton polygon} of this
tropical curve.

Let $K$ be the field of Puiseux series with complex coefficients
equipped with its standard non-Archimedian valuation $\val: K^*
\to \R$. According to Kapranov's theorem~\cite{Ka}, the underlying
set $\Pi$ of the tropical curve associated with $(A, \nu)$ is the
closure in $\R^2$ of the {\it non-Archimedian amoeba}
(see~\cite{Ka} and~\cite{Mi1}) of a curve in $(K^*)^2$ defined by
a polynomial $f(z,w)=\sum_{(i,j)\in A}c_{i,j}z^iw^j$ such that for
any point $(i,j)$ in $A$ one has $\val(c_{i,j})=\nu(i,j)$.
Clearly, such a polynomial~$f$ determines the pair $(A, \nu)$, and
thus, the associated tropical curve.

The function~$\nu$ defines a subdivision of $\Delta(A)$ into
convex polygons in the following way. Consider the overgraph
$\Gamma_\nu$ of~$\nu$, {\it i.e.}, the convex hull of the set
$\{(i,j,k) \in \R^3 : \; (i,j) \in A \;, \; k \geq \nu(i,j)\}$.
The polyhedron $\Gamma_\nu$ is naturally projected onto
$\Delta(A)$. The faces of $\Gamma_\nu$ which project injectively,
define a subdivision of $\Delta(A)$. Denote this subdivision by
$S(A, \nu)$. Let~$T$ be the tropical curve associated with $(A,
\nu)$. Note that $T$ does not determine uniquely the pair $(A,
\nu)$ (and even the polygon $\Delta(A)$). However, once the
polygon $\Delta(A)$ is fixed, the tropical curve~$T$ determines
uniquely the subdivision $S(A,\nu)$. For any tropical curve~$T$
with Newton polygon $\Delta$, denote by ${\cal D}_{\Delta}(T)$ the
subdivision of $\Delta$ determined by $T$.

A tropical curve is called {\it irreducible} if it cannot be
presented as a union of two proper tropical subcurves. A tropical
curve~$T$ with Newton polygon $\Delta$ is called {\it nodal} if
the subdivision ${\cal D}_\Delta(T)$ of $\Delta$ verifies the
following properties.
\begin{itemize}
\item any polygon of
${\cal D}_\Delta(T)$ is either triangle or a parallelogram,
\item any integer point on the boundary of
$\Delta$ is a vertex of ${\cal D}_\Delta(T)$.
\end{itemize}
In this case, the subdivision ${\cal D}_\Delta(T)$ is also called
nodal. Let $T$ be a nodal tropical curve with Newton polygon
$\Delta$. The {\it rank} of~$T$ is the difference diminished by
$1$ between the number of vertices of ${\cal D}_\Delta(T)$ and the
number of parallelograms in ${\cal D}_\Delta(T)$. The {\it
multiplicity} $\mu(T)$ of $T$ (and the multiplicity $\mu({\cal
D}_\Delta(T))$ of ${\cal D}_\Delta(T)$) is the product of areas of
all the triangles in ${\cal D}_\Delta(T)$ (we normalize the area
in such a way that the area of a triangle whose only integer
points are the vertices is equal to $1$).

Let $n$ be a natural number, $\bu$ a generic set of $n$ points in
$\R^2$. Consider the collection ${\cal C}(\bu)$ of nodal tropical
curves of rank $n$ and of Newton polygon $\Delta$ which pass
through all the points of $\bu$, and denote by ${\cal
T}^\C_n(\bu)$ the number of curves in ${\cal C}(\bu)$ counted with
their multiplicities.

\begin{theorem}\label{complex}{\rm (G.~Mikhalkin, see~\cite{Mi}
and~\cite{Sh})}. Let $\bu$ be a generic set of $n=r(\Delta)+g$
points in $\R^2$, where $0 \leq g \leq \delta(\Delta)$ is an
integer. Then ${\cal T}^\C_{n}(\bu)$ is equal to the number of
complex curves in the linear system $|D|$ which pass through a
fixed generic collection of $n$ points in $\Sig$ and have
$\delta(\Delta)-g$ nodes.
\end{theorem}

A nodal tropical curve~$T$ with Newton polygon $\Delta$ and the
corresponding subdivision ${\cal D}_\Delta(T)$ of $\Delta$ are
called {\it odd}, if each triangle in ${\cal D}_\Delta(T)$ has an
odd area. A nodal tropical curve~$T$ with Newton polygon $\Delta$
and the subdivision ${\cal D}_\Delta(T)$ are called {\it positive}
(resp., {\it negative}) if the sum of the numbers of interior
integer points over all the triangles of ${\cal D}_\Delta(T)$ is
even (resp, odd). Denote by ${\cal T}^+_n(\bu)$ (resp., ${\cal
T}^-_n(\bu)$) the number of odd positive (resp., negative) {\bf
irreducible} curves in ${\cal C}(\bu)$ (counted without
multiplicities).

\begin{theorem}\label{Welsch}{\rm (cf.~\cite{Mi} and~\cite{Sh})}.
Let $\bu$ be a generic set of $n=r(\Delta)$ points in $\R^2$. Then
$W_{\Sigma,D} = {\cal T}^+_{n}(\bu) - {\cal T}^-_{n}(\bu)$.
\end{theorem}

Theorem~\ref{Welsch} immediately follows from Proposition 5.1 and
Lemma 2.6 in~\cite{Sh}.

\medskip

{\bf C. Lattice paths.} Denote by~$p$ the highest
point of $\Del$ on the vertical axis,
and by~$q$
the most right point of $\Del$ on the horizontal axis.
The points $p$ and $q$ divide the boundary of $\Delta$ in two
parts. Denote the upper part of the boundary by $\partial
\Delta_+$, and the lower one by $\partial \Delta_-$. Fix a linear
function $\lambda: \R^2 \to \R$ which is injective on the integer
points of $\Delta$, and such that the restriction of $\lambda$ on
$\Delta$ takes its minimum
at~$p$ and takes
its maximum
at~$q$.

Let $l$ be a natural number. A path $\gamma:[0,l]\to\Delta$ is
called {\it $\lambda$-admissible} if
\begin{itemize}
\item $\gamma(0) = p$ and $\gamma(l) = q$,
\item the composition $\lambda \circ \gamma$ is injective,
\item for any integer $0 \leq i \leq l-1$ the point $\gamma(i)$
is integer, and $\gamma([i,i+1])$ is a segment.
\end{itemize}
The number~$l$ is called the {\it length} of~$\gamma$. A
$\lambda$-admissible path~$\gamma$ divides $\Delta$ in two parts:
the part $\Delta_+(\gamma)$ bounded by $\gamma$ and $\partial
\Delta_+$ and the part $\Delta_-(\gamma)$ bounded by $\gamma$ and
$\partial \Delta_-$. Define an operation of {\it compression} of
$\Delta_+(\gamma)$ in the following way. Let $j$ be the smallest
positive integer $1 \leq j \leq l - 1$ such that $\gamma(j)$ is
the vertex of $\Delta_+(\gamma)$ with the angle less than $\pi$ (a
compression of $\Delta_+(\gamma)$ is defined only if such an
integer~$j$ does exist). A compression of $\Delta_+(\gamma)$ is
$\Delta_+(\gamma')$, where $\gamma'$ is either the path defined by
$\gamma'(i)=\gamma(i)$ for $i < j$ and $\gamma'(i)=\gamma(i+1)$
for $i \geq j$, or the path defined by $\gamma'(i)=\gamma(i)$ for
$i \ne j$ and $\gamma'(j)=\gamma(j-1)+\gamma(j+1)-\gamma(j)$ (if
$\gamma(j-1)+\gamma(j-1)+\gamma(j) \in \Delta$). Note that
$\gamma'$ is also a $\lambda$-admissible path. A sequence of
compressions started with $\Delta_+(\gamma)$ and ended with a path
whose image coincides with $\partial \Delta_+$ defines a
subdivision of $\Delta_+(\gamma)$ which is called {\it
compressing}. A compression and a compressing subdivision of
$\Delta_-(\gamma)$ is defined in the completely similar way. A
pair $({\cal S_+(\gamma)}, {\cal S_-(\gamma)})$, where ${\cal
S_\pm(\gamma)}$ is a compressing subdivision of
$\Delta_\pm(\gamma)$, produces a subdivision of $\Delta$. Denote
by ${\cal N}_\lambda(\gamma)$ the collection of nodal subdivisions
of $\Delta$ which can be obtained in this way starting
with~$\gamma$.

\begin{theorem}\label{paths}{\rm (G.~Mikhalkin, see~\cite{Mi})}.
Let $0 \leq g \leq \delta(\Delta)$ be an integer. There exists a
generic set $\bu$ of $n=r(\Delta) + g$ points in $\R^2$ such that
the map ${\cal D}_\Delta$, associating to a nodal tropical curve
$T$ with Newton polygon $\Delta$ the corresponding subdivision
${\cal D}_\Delta(T)$ of $\Delta$ {\rm (}see subsection B{\rm )},
establishes a 1-to-1 correspondence between the set ${\cal
C}(\bu)$ and the disjoint union $\amalg_\gamma {\cal
N}_\lambda(\gamma)$, where $\gamma$ runs over all the
$\lambda$-admissible paths in $\Delta$ of length~$n$. In
particular, ${\cal T}^\C_n(\bu)=\sum_\gamma \sum_{{\cal S}\in
{\cal N}_\lambda(\gamma)} \mu({\cal S})$, where $\mu({\cal S})$ is
the multiplicity of ${\cal S}$.
\end{theorem}

Take a set $\bu$ with the properties described in
Theorem~\ref{paths}, and denote by ${\cal N}^+_\lambda(\gamma)$
(resp., ${\cal N}^-_\lambda(\gamma)$) the number of odd positive
(resp., negative) subdivisions in ${\cal N}_\lambda(\gamma)$ which
are the images of {\bf irreducible} tropical curves under the
bijection ${\cal D}_\Delta\mid_{{\cal C}(\bu)}: {\cal C}(u) \to
\amalg_\gamma {\cal N}_\lambda(\gamma)$. The following statement
is an immediate corollary of Theorems~\ref{Welsch}
and~\ref{paths}.

\begin{proposition}\label{new-proposition}
The Welschinger number $W_{\Sigma, D}$ is equal to
$\sum_\gamma({\cal N}^+_\lambda(\gamma) - {\cal
N}^-_\lambda(\gamma))$, where $\gamma$ runs over all the
$\lambda$-admissible paths in $\Delta$ of length~$r(\Delta)$.
\end{proposition}

\medskip

{\bf D. Positivity of the Welschinger invariant.} Let $\lambda^0:
\R^2 \to \R$ be a linear function defined by
$\lambda^0(i,j) = i - \varepsilon j$,
where $\varepsilon$ is
sufficiently small positive irrational number (so that $\lambda^0$
defines a kind of a lexicographical order on the
integer
points
of the polygon $\Delta$).

\begin{proposition}\label{p2}
For any $\lambda^0$-admissible path $\gamma$ in $\Delta$, the
number ${\cal N}^-_{\lambda^0}(\gamma)$ equals $0$.
\end{proposition}

{\bf Proof.} Let $\gamma$ be a $\lambda^0$-admissible path in
$\Delta$, and ${\cal S}$ a subdivision in the collection
${\cal N}_{\lambda^0}(\gamma)$.
The subdivision ${\cal S}$ does not have an edge with the endpoints
$(i_1, j_1)$ and $(i_2, j_2)$ such that $|i_1 - i_2|> 1$;
otherwise, at least one integer point on the boundary of~$\Delta$
would not be a vertex of the corresponding compressing
subdivision. This implies that no triangle in ${\cal S}$ has
interior integer points.
\proofend

{\bf Proof of Theorems~\ref{t3} and~\ref{t3n}.}
We consider only the cases of $\C P^2$ and $P_3$, since the
constructions presented below can be easily adapted to the
remaining cases.

If $\Sigma = P_3$, we can assume that
$d_1+d_2+d_3\le d$. Indeed,
we have
$|dL-d_1E_1-d_2E_2-d_3E_3|=|d'L-d'_1E'_1-d'_2E'_2-d'_3E'_3|$,
where $d'=2d-d_1-d_2-d_3$, $d'_i=d-d_j-d_k$, $E'_i=L-E_j-E_k$ for
all $\{i,j,k\}=\{1,2,3\}$.
It remains to notice that
$d'_1+d'_2+d'_3<d'$ as far as $d_1 + d_2 + d_3 > 0$,
and $(d_1!)/(d-d_2-d_3)! = (d'-d'_2-d'_3)!/(d'_1!)$.

\begin{figure}
\begin{center}
\epsfxsize 114mm
\epsfbox{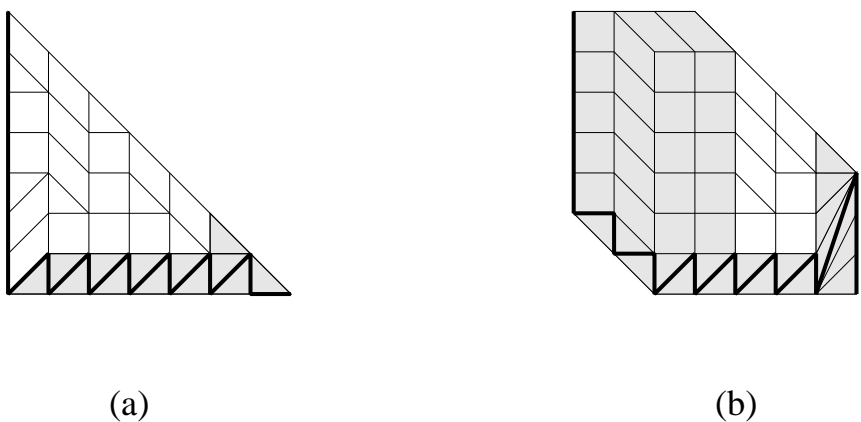}
\end{center}
\caption{The chosen $\lam^0$-admissible paths and examples of
subdivisions} \label{f1}
\end{figure}

Take the $\lambda^0$-admissible lattice path $\gamma$ of length
$r(\Delta)$ passing successively through
$r(\Del)+1$ integer points of
$\Del$ described below.
\begin{itemize}
\item If $\Sig=\C P^2$, then the points are
$(0,d-j)$, where $j=0,...,d$; \hskip10pt
$(i,1)$, $(i,0)$, where $i=1,...,d-1$; and
$(d,0)$ (see Figure \ref{f1}(a)).
\item If $\Sig = P_3$, then the points are $(0,d-d_2-j)$,
where $j=0,...,d-d_2-d_3$; \hskip10pt
$(i,d_3-i+1)$, $(i,d_3-i)$, where $i=1,...,d_3$; \hskip15pt
$(i,1)$, $(i,0)$, \hskip10pt where $i=d_3+1,...,d-d_1-1$; \hskip10pt
$(d-d_1,d_1-j)$,
where $j=0,...,d_1$
(see Figure~\ref{f1}(b)).
\end{itemize}

Then we choose
the following subdivisions of $\Del$ constructed along the procedure of
subsection C starting with $\gamma$.

If $\Sig=\C P^2$, in the domain
$\{(\xi,\eta)\in\Del\ :\
0\le\xi\le1,\ \eta\ge\xi\}$,
to which we attribute the number
$i=0$, and in each strip $\{(\xi,\eta)\in\Del\ :\ i\le\xi\le i+1,
\ \eta \ge 1\}$, $1\le i \le d-3$, we pack arbitrarily $d-i-2$
parallelograms of (normalized) area $2$.
The remaining part
of~$\Delta$ is uniquely covered by triangles of area~$1$
in order to obtain a pair $({\cal S_+(\gamma)}, {\cal S_-(\gamma)})$
of compressing subdivisions.

If
$\Sig=P_3$, in the domain
$$\big\{(\xi,\eta)\in\Del\ :\ \xi\le d_2+d_3-1,\ \eta\ge
\min\{d_3,\max\{d_3+1-\xi,1\}\}\big\}$$
we pack parallelograms of area~$2$
in a way compatible with the compressing procedure
(see Figure~\ref{f1}(b)).
Then in each strip $$\{(\xi,\eta)\in\Del\
:\ d_2+d_3-1+i\le\xi\le d_2+d_3+i,\ \eta\ge 1\}\ ,$$
$i=0,...,d-d_1-d_2-d_3-1$, we pack arbitrarily $d-d_2-d_3-i-1$
parallelograms of area~$2$. The remaining part of $\Del$ is uniquely
covered by triangles of area~$1$
in order to obtain a pair $({\cal S_+(\gamma)}, {\cal S_-(\gamma)})$
of compressing subdivisions.

We illustrate the construction
in Figure~\ref{f1}, where the shadowed domain
represents the common part of the described subdivisions.

All the chosen subdivisions
of $\Del$ are nodal and odd, and they correspond to
irreducible tropical curves. The number of subdivisions under
consideration is equal respectively to $d!/2$ for $\C P^2$,
and to $(d-d_2-d_3)!/d_1!$ for $P_3$.
Thus, according to
Propositions~\ref{new-proposition} and~\ref{p2} the same numbers
are lower bounds for the Welschinger invariant. \proofend

\section{Non-invariance of the Welschinger number for elliptic
curves}\label{sec4}

\begin{theorem}\label{t2} For any $d\ge 4$,
there exist generic collections $\bw'$ and $\bw''$ of $3d$
distinct points in $\R P^2$ such that
$W_{\C P^2,d,1}( \bw')\ne W_{\C P^2,d,1}(\bw'')$.
\end{theorem}

{\bf Proof}. As is known (see, for example, \cite{Se}),
for any generic set $\bz$ of $3d-1$ points in
$\C P^2$, the plane nodal irreducible elliptic curves of degree $d$,
which pass through each of the points $z_1,\dots,z_{3d-1}$ in
$\bz$ and are nonsingular at these points, form a smooth
quasi-projective one-dimensional variety
$V=V( \bz)\subset|{\cal O}_{\C P^2}(d)|\simeq \C P^{d(d+3)/2}$;
the Zariski (projective)
tangent space to $V$ at $C\in V$ is the one-dimensional linear
system $\Lam(\Sing(C),\bz)$ of curves of degree $d$ which pass
through
$\Sing(C)$ and
$\bz$.

Denote by $\zeta:\tilde C\to C$ the normalization of $C\in V$, by
$\tilde \bz=\{\tilde z_1,\dots, \tilde z_{3d-1}\}$ and $\tilde S$
the pull-back of $\bz$ and $\Sing (C)$, and by ${\cal N}$ the
normal line bundle of the immersion $\displaystyle{\phi:\tilde C
\mathto^\zeta C\to \C P^2}$. According to the adjunction formula,
the bundle ${\cal N}$ is of degree $3d$ and is induced by
${\cal O}_{\C P^2}(3)$, and, as follows from Riemann-Roch formula, there is
one and only one point $\tilde\omega_{C,\bz}\in\tilde C$ such that
$[\tilde z_1+\dots+\tilde z_{3d-1}+\tilde
\omega_{C,\bz}]=c_1({\cal N})\in\Pic \ \tilde C$. On the other
hand, for any $C^{(1)}\ne C$ from $\Lam(\Sing(C),\bz)$ one has
\begin{equation}
\phi^* C^{(1)}=\tilde z_1+\dots+\tilde z_{3d-1}+\tilde
\omega_{C,\bz}+\sum_{s\in \tilde S} s,
\label{e1n}\end{equation} and any such curve $C^{(1)}$ of degree
$d$ together with $C$ generates $\Lam(\Sing(C),\bz)$.

If $z_1,\dots,z_{3d-1}$ are real, the variety $V$ is real, and for
any real $C\in V$ all the above objects,
including $\tilde C$
and
$\tilde\omega_{C,\bz}$, are real as well.

\medskip

{\bf Claim A}.
For a generic $\bz$
the map
$C\in V\mapsto
\omega_C=\phi(\tilde\omega_{C,\bz})\in \C P^2$
is not constant.

\medskip

Denote by $\bar V$ the closure of $V$, by $\bar{\cal V}$ the
closure of ${\cal V} =\{(C, p) \ \in V\times \C P^2 \ :\ \; \; p\in C\}$,
and by $\pi: \bar V\times \C P^2 \to \C P^2$ the canonical
projection. Here, we deal with an immersed surface $\cal V$, and,
according to the above characterization of
$\Lambda(\Sing(C),\bz)$,
the image $\theta$ of the map $C\in V\mapsto \omega_C\in \C P^2$
is the critical value set of
$\pi\big|_{\cal V}$,
whereas $\Theta=\{(C, \w_{C,\bz})\ :\ C\in V\}\subset{\cal V}$
is its critical point set.

\medskip

{\bf Claim B}. For a generic $\bz$
and a generic $\omega_C\in\theta$, the fibre $\pi^{-1}(\omega_C)$
meets $\bar{\cal V}$ only inside ${\cal V}$, and, in addition, all
the points $(C',\omega_C)\in\pi^{-1}(\omega_C)$ distinct from
$(C,\omega_C)$ are regular.

\medskip

Hence the set $\pi^{-1}(z)\cap{\cal V}$ may bifurcate only in a
neighborhood of $(C,\omega_C)$ as $z$ ranges over a neighborhood
of $\omega_C$ in $\R P^2$.

Let $\w$ be a smooth point of $\theta$ and let $\omega=\omega_{C}$
for some $C\in V$. The order $m$ of ramification of $\pi$ at
$(C,\w)$ is given by the variation formula
$$
m=m(\tilde C(t)\cdot \Theta(t))_{\tilde\w(C,\bz)}=(C(t)\cdot\theta
(t))_{\w(C(t))}= 1+(\frac{d}{dt}C(t)\cdot C(t))_{\w(C(t))},
$$
where $t$ is a uniformization parameter of $\theta$ at $\w$. On
the other hand,
B\'ezout's theorem implies ({\it cf}. (\ref{e1n})) that
$(\frac{d}{dt}C(t)\cdot C(t))_{\w(C(t))}=1$.
Since $\pi\big|_{\cal V}$ has a double folding
along the germ of $\theta$ at $\omega$, the set
$\pi^{-1}(z)\cap\R{\cal V}$ (where $\R{\cal V}$ is the real part
of $\cal V$) either looses or obtains two elements when $z$ varies
in $\R P^2$ along a line through $\omega_C$ transverse to the real
part $\R C$ of $C$. Both the elements of
$\pi^{-1}(z)\cap\R{\cal V}$ which appear (disappear) in this
bifurcation, are real nodal elliptic curves with the same number
of solitary nodes as $C$, and hence the Welschinger number
changes. Theorem \ref{t2} is proven. (The proofs of Claims A and B
are found in Appendix. ) \proofend

\begin{remark}\label{r2} The above bifurcation can be
characterized as a jump of $\dim\Lam(\Sing(C),z_1,...,z_{3d})$
from $0$ to $1$, where $\Lam(\Sing(C),z_1,...,z_{3d})$ is the
linear system of curves of degree $d$, passing through
the singular locus $\Sing(S)$
of a nodal elliptic curve $C$ of degree $d$ and through some
points $z_1,...,z_{3d}\in C\backslash\Sing(C)$. This cannot happen
in the case of rational curves, since, by Riemann-Roch theorem,
$\dim\Lam(\Sing(C),z_1,...,z_{3d-1})=0$ for any rational nodal
curve $C$ of degree $d$ and any distinct points
$z_1,...,z_{3d-1}\in C\backslash\Sing(C)$.
\end{remark}

\section{A remark on counting real plane curves with one node}
\label{sec5}

To count curves with one node, one has to consider linear pencils.
In the case of plane curves of degree $d$, the computation of the
Euler characteristic of $\R P^2$ by means of a pencil with $d^2$
real base points as done in \cite{DK} for $d=3$, gives $d^2-1$ as
the lower bound for the number of real nodal curves in such a
pencil. The upper bound $3(d-1)^2$ is given by the number of
complex nodal curves in the pencil. To show that the bounds are
sharp, consider a pencil generated by
$\prod_{\alpha=1}^d(x-x_\alpha)$ and
$\prod_{\alpha=1}^d(y-y_\alpha)$. The number of real nodal curves
in this pencil is $(d-1)^2$. If we perturb
$\prod_{\alpha=1}^d(x-x_\alpha)$ and
$\prod_{\alpha=1}^d(y-y_\alpha)$ into unions of generic lines and
then take a close generic pencil, the number of real nodal curves
becomes $3(d-1)^2$. If we perturb $\prod_{\alpha=1}^d(x-x_\alpha)$
and $\prod_{\alpha=1}^d(y-y_\alpha)$ by means of affine
polynomials without minima and maxima (see \cite{Shustin}), we
obtain $(d-1)^2+2(d-1)=d^2-1$
real nodal curves in the pencil.

Welschinger noticed that his number is not invariant for curves
with one node if $d\ge 4$. The argument of Welschinger is the
following. A pencil of plane curves of degree $d$ is defined by
$\frac12 d(d+3) -1$ fixed points, and it has $\frac{1}{2}d(d-3)+1$
extra base points. From the above computation of the Euler
characteristic, it follows that, in this case, the Welschinger
number is equal to $d^2-1-2r$, where $r$ is the number of pairs of
conjugated imaginary points among the base ones. It remains to
consider pencils with $r$ varying from $0$ to $[d(d-3)/4+1/2]$ and
to choose each time all the defining $\frac12 d(d+3) -1$ points
real.

\bigskip

{\bf Acknowledgements}. The first and the second authors are
members of Research Training Network RAAG CT-2001-00271;
the second author is as well a
member of Research Training Network EDGE CT-2000-00101; the third author has
been supported by the Bessel research award from the Alexander von
Humboldt Stiftung/Foundation. The third author also thanks
Universit\'e Louis Pasteur, Strasbourg, Universit\'e de Rennes I
and Universit\"at Kaiserslautern for the hospitality and excellent
work conditions.

\bigskip

\section*{Appendix}\label{ap}

{\it Proof of Claim A}. Since the non-constancy of $\omega_C$ is
an open condition, it is sufficient to verify it for a suitable
special configuration of points. Let us choose ten generic points
$z_1,...,z_{10}$, trace a (unique) cubic $C_3$ through
$z_1,...,z_9$, and pick a rational nodal curve $R_{d-3}$ of degree
$d-3$, passing through $z_{10}$ and intersecting transversely
$C_3$ at $3(d-3)$ points, which all are nonsingular for $R_{d-3}$.
Denote by
$z'_{11},...,z'_{3d-1},p_1,p_2$ the latter intersection points.
Consider the elements $C\in
V=V(z_1,\dots,z_{10},z'_{11},...,z'_{3d-1})$ which are close to
$C_3R_{d-3}$, which have nodes near $p_1$, near each of $z'_j$
($j=11,...,3d-1$), and near $\Sing(R_{d-3})$, and for which the
branch passing through $z'_j$ ($j=11,...,3d-1$) approaches
$R_{d-3}$ as $C$ tends to $C_3R_{d-3}$. By the standard
transversality arguments, they smooth off the node at $p_2$ and
form, together with $C_3R_{d-3}$ a smooth (germ of)
one-dimensional variety $\widetilde V$ (one can easily verify the
first order conditions which ensure that the linear system
$\Lam(p_2)$ of curves of degree $d$ through $p_2$, the linear
systems of curves of degree $d$ meeting $R_{d-3}$ at $z'_j$ with
multiplicity $\ge 2$ and the linear system
$\Lam(z_1,...,z_{10},\Sing(R_{d-3}))$ intersect transversely all
together, cf. \cite{Se} and \cite{Sh0}).

From now on we argue to contradiction. Since nodal elliptic curves
of given degree form an irreducible variety (see \cite{Har}), for
Claim A it suffices to treat the case when $\omega_C$ is a
constant point $\w$ for all $C\in\widetilde V
\backslash\{C_3R_{d-3}\}$. Then $\w$ is the intersection point of
$C_3R_{d-3}\backslash(\{z_1,...,z_{10},z'_{11},...,z'_{3d-1},p_1\}
\cup\Sing(R_{d-3}))$ with a curve $C'=\frac{d}{dt}C(t)\in
\Lam(z_1,...,z_{10},z'_{11},...,z'_{3d-1},p_1)
\backslash\{C_3R_{d-3}\}$, where $C(t)$ is a parametrization of
$\widetilde V$. On the other hand, $C\in\widetilde V
\backslash\{C_3R_{d-3}\}$ meets $C_3$ at
$z_1,...,z_9,z'_{11},...,z'_{3d-1},w$ and some point $p'_1$ in a
neighborhood of $p_1$. Thus,
$$[z_1+...+z_9+z'_{11}+...+z'_{3d-1}+p'_1+w]=
[z_1+...+z_9+z'_{11}+...+z'_{3d-1}+p_1+w]\in \Pic(C_3)\ ,$$ and it
remains to check that $p'_1\ne p_1$ under suitable choice of
$z_1,\dots,z_9$.

Suppose further that all $C\in\widetilde V$
pass through $p_1$. Then $C'=\frac{d}{dt}C(t)$ belongs either to
the space $T'$ of curves of degree $d$ crossing $C_3$ at $p_1$
with multiplicity $\ge 2$, or to the space $T''$ of curves of
degree $d$ crossing $R_{d-3}$ at $p_1$ with multiplicity $\ge 2$.
However, the intersection of $T'$ or $T''$ with the linear system
$\Lam'$ of curves of degree $d$ passing through
$z_1,...,z_{10},\Sing(R_{d-3})$ and crossing $R_{d-3}$ at each of
$z'_{11},...,z'_{3d-1}$ with multiplicity $\ge 2$, is transversal
and zero-dimensional, since contains only $C_3R_{d-3}$. For
example, a curve $H$ of degree $d$ from $T'\cap\Lam'$ crosses
$C_3$ at $z_1,...,z_9,z'_{11},...,z'_{3d-1}$ and twice at $p_1$,
but
$$z_1+...+z_9+z'_{11}+...+z'_{3d-1}+2p_1$$
can not be cut on $C_3$ by a curve of degree $d$, at least if we
move slightly one of the points $z_1,\dots,z_9$ along $C_3$. Hence
$H$ splits off $C_3$. In turn the remaining component of degree
$d-3$ coincides with $R_{d-3}$, since it meets $R_{d-3}$ at
$z_{10},z'_{11},...,z_{3d-1},\Sing(R_{d-3})$ with the total
multiplicity $\ge 3d-10+(d-4)(d-5)=(d-3)^2+1.$\proofend

\medskip

{\it Proof of Claim B}. To prove the first part of Claim B, it is
sufficient to show that for a generic $\bz$ no curve
$C'\in\overline{V(\bz)}\backslash V(\bz)$ has a component in
common with $\theta=\theta(\bz)$. So assume by contrary that, for
a generic $\bz$, some $C'\in\overline{V(\bz)}\backslash V(\bz)$
has a component $C'_1\subset\theta$. Pick a generic point $q\in
C'_1$ and a curve $C\in V(\bz)$ with $\omega_C= q\ne z_1$ (recall
that the map $C\in V\mapsto \omega_C$ is not constant). Due to the
generality of $\bz$,
we may suppose that at least for small variations of $\bz$
the variations of $\overline{\cal V}\setminus\cal V$ are flat.
Consider a deformation $z_1(t)$, $t\in(\C,0)$, of $z_1$ along
$C'$. Due to the above flatness, $C'_1$ remains a component of the
curve $\theta(\bz(t))$, $\bz(t)=\{z_1(t),z_2,...,z_{3d-1}\}$, and
$C=C(0)$ deforms in a family $C(t)$ of elliptic curves passing
through $\bz(t)$ and satisfying
$\omega_{C(t)}=q$. By construction, $\theta$ is an envelope of the
family $V$, and hence $C'_1$ is tangent to each of $C(t)$,
$t\in(\C,0)$, at $q$. Thus,
\begin{equation}C(t)\cap C
= z_2+...+z_{3d-1}+ 2q+D(t)
\end{equation} with a divisor $ D(t)$ close to $ D=\phi(\tilde
S)$. The limit $t=0$ yields the relation \begin{equation} [\tilde
z_2+...+\tilde z_{3d-1}+2\tilde q
+
\tilde S]=
\phi^*
(O_{P^2}(d)) \in\Pic(\tilde C)\ . \label{e2n}
\end{equation}
On the other hand, (\ref{e1n}) implies
\begin{equation}[\tilde z_1+...+\tilde z_{3d-1}+\tilde w_{C}+
\tilde S]=
\phi^*
(O_{P^2}(d)) \in\Pic(\tilde C)\ ,
\end{equation}
which contradicts (\ref{e2n}) since $q=\omega_C$ and $[z_1]\ne
[\omega_C]\in\Pic(C)$.

We prove the regularity statement arguing to
contradiction
as
well. Since the nodal elliptic curves of given degree form an
irreducible variety (see
\cite{Har}),
it suffices to treat
the case when for generic $\bz$ the whole curve
$\Theta\subset{\cal V}$ is multiply mapped on $\theta$ by $\pi$.
Then, for a generic $C\in V(\bz)$ there exists $\hat C\in V(\bz),
\hat C\ne C,$
with a common point $\omega_C=\omega_{\hat C}$ on a nonsingular
branch of $\theta$. Both $C$ and $\hat C$ are
tangent to
$\theta$ at this point and
such a
configuration must
be
preserved under small
variations of $\bz$.

Note that $C$ and $\hat C$ intersect transversally at least at one
of the points $z_1,...,z_{3d-1}$. Indeed, otherwise, taking a
deformation $z_1(t),z_2(t)\in C$, $t\in(\C,0)$, such that
$[z_1(t)+z_2(t)]=[z_1+z_2]\in\Pic(C)$, we obtain that $\omega_C$
does not depend on $t$, and $\hat C=\hat C(0)$ varies in a family
$\hat C(t)$ of curves passing through
$z_1(t),z_2(t),z_3,...,z_{3d-1}$ and having $\omega_{\hat
C(t)}=\omega_C$. However, the tangency assumption would imply that
the total intersection of $\hat C$ and $\hat C(t)$, $t\ne 0$, at
$z_3,...,z_{3d-1},\omega_C$ and in a neighborhood of $\Sing(C)$ is
at least $6d-4+d(d-3)=d^2+3d-4>d^2$, which contradicts B\'ezout's
theorem.

Thus, we can suppose that $C$ and $\hat C$ are
nonsingular and transverse at $z_1$. The
first order infinitesimal equisingular deformations of the
curve
$\hat C$ with fixed points $z_2,...,z_{3d-1}$ correspond
to the elements of the
two-dimensional space $H^0(\hat C,{\cal J}_{Z/\hat C}\otimes{\cal
O}_{\C P^2}(d))$, where ${\cal J}_{Z/\hat C}$ is the ideal sheaf of
the zero-dimensional scheme $Z\subset \hat C$ defined by the
maximal ideals at $z_2,...,z_{3d-1}$ and $\Sing(\hat C)$. Pick
a section $s\in H^0(\hat C,{\cal
J}_{Z/\hat C}\otimes{\cal O}_{\C P^2}(d))$ which does not vanish
neither at $z_1$ nor at $\omega_{\hat C}$. Consider a deformation
$\hat C(t)$, $t\in(\C,0)$, $\hat C(0)=\hat C$, generated by $s$,
{\it i.e.}, $\frac{d}{dt}\hat C(t)\big|_0=s$, and keep $C$ fixed,
{\it i.e.}, put $C(t)=C$.
The transversality of $C$ and $\hat C$ at $z_1$ implies that, for
small $t$, there is a unique intersection point $z_1(t)$ of $\hat
C(t)$ and $C$ near $z_1$.
Put $\bz(t)=\{z_1(t),z_2,\dots,z_{3d-1}\}$.
Since $s(\omega_C)\ne 0$,
we
have $\frac{d}{dt}\omega_{\hat C(t), \bz(t)}\notin T\theta$ (where
$T$
stands for the tangent space). On the other hand,
$\frac{d}{dt}\omega_{C(t),\bz(t)}\in T\theta$.

It means that the branches of $\theta(\bz)$ at
$\omega_C=\omega_{\hat C}$ split into separate branches in the
above variation, which contradicts the initial assumption on the
stability of the multiple projection $\Theta(\bz)\to\theta(\bz)$.
\proofend

\bigskip
\vskip40pt

{\ncsc CNRS \\[-21pt]

Institut de Recherche Math\'ematique de Rennes \\[-21pt]

Campus de Beaulieu, 35042 Rennes Cedex, France} \\[-21pt]

{\it E-mail address}: {\ntt Ilia.Itenberg@univ-rennes1.fr}

\vskip10pt

{\ncsc Universit\'e Louis Pasteur et IRMA \\[-21pt]

7, rue Ren\'e Descartes, 67084 Strasbourg Cedex, France} \\[-21pt]

{\it E-mail address}: {\ntt kharlam@math.u-strasbg.fr}

\vskip10pt

{\ncsc School of Mathematical Sciences \\[-21pt]

Raymond and Beverly Sackler Faculty of Exact Sciences\\[-21pt]

Tel Aviv University \\[-21pt]

Ramat Aviv, 69978 Tel Aviv, Israel} \\[-21pt]

{\it E-mail address}: {\ntt shustin@post.tau.ac.il}

\end{document}